\documentclass[11pt,a4paper,twoside]{amsart}
\usepackage[T1]{fontenc}
\usepackage{mathtools,microtype}
\usepackage{txfonts}
\usepackage[utf8]{inputenc}
\usepackage{amsthm}
\usepackage{amsmath}
\usepackage{amssymb}
\usepackage{enumitem}
\usepackage{xcolor}
\usepackage{comment}
\usepackage{tikz}
\usepackage[colorlinks,citecolor=blue,linkcolor=blue,urlcolor=blue]{hyperref}

\makeatletter
\renewcommand{\subsection}{\@startsection
  {subsection}{2}{\z@}%
  {-3.25ex\@plus -1ex \@minus -.2ex}%
  {1.2ex \@plus .2ex}%
  {\normalfont\bfseries}}
\makeatother

\makeatletter
\renewcommand{\subsubsection}{\@startsection
  {subsubsection}{3}{\z@}%
  {-3.25ex\@plus -1ex \@minus -.2ex}
  {1.2ex \@plus .2ex}
  {\normalfont\bfseries}}
\makeatother

\usetikzlibrary{arrows.meta}

\theoremstyle{plain}
\newtheorem{theorem}{\bfseries Theorem}[section]
\newtheorem{lemma}[theorem]{\bfseries Lemma}

\newtheorem*{introtheorem}{\bfseries Theorem}
\newtheorem*{conjecture}{\bfseries Conjecture}

\theoremstyle{definition}
\newtheorem{definition}[theorem]{\bfseries Definition}

\newtheorem{remark}[theorem]{\bfseries Remark}

\newcommand{\rint}[0]{\mathcal{O}_F}
\newcommand{\disc}[0]{D_F}
\newcommand{\facets}[0]{\mathcal{F}^1(\sigma)}
\newcommand{\kfacets}[0]{\mathcal{F}^k(\sigma)}

\newcommand{\perfect}[0]{\mathcal{P}_}
\newcommand{\GL}[0]{\operatorname{GL}_n(\mathbb{Z})}

\newcommand{\sym}[0]{\operatorname{Sym}_n(\mathbb{R})}
\newcommand{\GLo}[0]{\operatorname{GL}_n({\mathcal{O}_F})}

\newcommand{\Her}[0]{\mathcal{H}_n(\fr)}

\newcommand{\vol}[0]{\mathrm{Vol}}

\newcommand{\fr}[0]{F_\mathbb{R}}
\newcommand{\Diff}{\mathfrak D_F}
\newcommand{\Codiff}{\mathfrak D_F^{-1}}
\newcommand{\tr}{\operatorname{Tr}}

\renewcommand{\le}{\leqslant}
\renewcommand{\ge}{\geqslant}
\renewcommand{\leq}{\leqslant}
\renewcommand{\geq}{\geqslant}

\title[{\SMALL Perfect Forms and Torsion Bounds in Algebraic $K$)-Theory }]{Perfect Forms and Torsion Bounds in Algebraic $K$-Theory over Totally Real and CM Fields}

\author{Alejandro de la Torre Durán}

\address{Univ. Grenoble Alpes, CNRS, Institut Fourier, F-38000 Grenoble, France}
\email{alejandro.de-la-torre-duran@univ-grenoble-alpes.fr}

\subjclass[2020]{11H55, 11F75, 11F06, 11Y99, 19D50, 20J06, 55N91}
\keywords{Perfect forms, Voronoi-Koecher complex, group cohomology, modular groups, Steinberg modules, K-theory of integers, well-rounded lattices, tessellations}

\begin{document}
\begin{abstract}
   We extend the method of Van Woerden \cite{woerden} to totally real and CM number fields to obtain an upper bound for the number of $\GLo-$equivalence classes of perfect forms. Moreover, by considering appropriate barycentric subdivisions of the Voronoi tessellation, we extend these estimates to the number of $\GLo$-orbits of $k$-dimensional faces of the Voronoi tessellation.
 As a consequence, we improve the asymptotic order of the bounds for the torsion part of $K_m(\rint)$ obtained by Soulé \cite{soule_hermitian_bounds} and Bayer-Fluckiger, Emery, and Houriet \cite{emery,Juliet_PHD}. We emphasise the improvement of the general bound for $F=\mathbb Q$, unchanged since Soulé's work, and the independence of the discriminant exponent from the field degree for totally real and CM fields.
\end{abstract}
\maketitle

\section*{Introduction}

The enumeration of perfect forms is a classical problem in the geometry of numbers, closely related to the lattice sphere-packing problem. Voronoi introduced an algorithm that solves the lattice sphere-packing problem in any fixed dimension by enumerating perfect forms \cite{Voronoi1908}. This naturally raises the question of how the number of equivalence classes of perfect forms grows with the dimension.

Let $p_{\mathbb{Q},n}$ denote the number of equivalence classes of perfect forms in dimension $n$, up to the action of $\GL$ and positive scaling. Soulé \cite{soule_bounds} established the upper bound
$$
p_{\mathbb{Q},n}\leq e^{O(n^4\log n)}.
$$
Bacher \cite{bacher} subsequently improved this estimate to $e^{O(n^3\log n)}$, while also establishing an exponential lower bound. Later, van Woerden \cite{woerden} used volumetric arguments to obtain the sharper upper bound
$$
p_{\mathbb{Q},n}\leq e^{O(n^2\log n)}.
$$
More recently, Moshchevitin \cite{moshchevitin2026} proved that
$$
p_{\mathbb Q,n} \geq e^{\left(\frac12 n^2\log n-O(n^2)\right)},
$$
showing that van Woerden's upper bound is asymptotically optimal in the exponent.

The work of Voronoi \cite{Voronoi1908} was later generalized by Koecher to arbitrary number fields \cite{Koecher}. Nevertheless, some properties of perfect forms in the rational case are lost. For instance, extremality does not imply perfection for $F\neq\mathbb Q$, and thus a solution to the lattice sphere-packing problem need not be among the perfect lattices. Moreover, much less is known about the enumeration of perfect forms over general number fields. While explicit classifications have been obtained for certain fields and small ranks, general estimates remain comparatively scarce. For instance, Scheckelhoff, Thalagoda and Yasaki \cite{scheckelhoff_thalagoda_yasaki} studied perfect binary Hermitian forms over imaginary quadratic fields, obtaining explicit classifications and a lower bound for the number of equivalence classes in terms of the discriminant of the field. 

A further motivation for studying perfect forms comes from Voronoi reduction theory and its applications to the topology of arithmetic groups. To each perfect form, one associates a full-dimensional polyhedral cone in the space $\Her$, generated by the rank-one forms corresponding to its minimal vectors. These cones, together with their faces, yield a polyhedral decomposition of a suitable rational closure of the cone of positive-definite forms. After projectivization, this gives an arithmetic-group-invariant CW decomposition \cite{Coulangeon}. This construction has been used to compute the cohomology of $\GL$ \cite{philippe_advances} and to prove the triviality of $K_8(\mathbb Z)$ \cite{k8}. Similar computations have also been carried out over arbitrary number fields \cite{hermitian, Vor_comp_num_fields}.

Beyond explicit cohomological computations, the geometry of the Voronoi decomposition can also be used to obtain quantitative bounds on torsion. In \cite{soule_bounds}, Soulé, by proving geometric and combinatorial estimates for the Voronoi tessellation, obtained bounds on torsion in the homology of arithmetic groups and, ultimately, in the algebraic $K$-groups of $\mathbb{Z}$.
This approach was subsequently extended by Soulé \cite{soule_hermitian_bounds} to the rings of integers of general number fields and later improved by Bayer-Fluckiger, Emery, and Houriet in \cite{emery}. In Houriet's doctoral thesis \cite{Juliet_PHD}, these bounds are stated with a slightly better asymptotic order, which is why we will use this reference when comparing our bounds with previous results. Nevertheless, the credit for these bounds belongs to all three aforementioned authors. More recently, Emery obtained improved non-explicit asymptotic bounds in terms of the discriminant for imaginary quadratic fields \cite{Emery_imaginary_quad}.

Fix a totally real or CM number field $F$ with ring of integers $\rint$, degree $d=[F:\mathbb Q]$, and discriminant $\disc$. The aim of this paper is to extend the quantitative study of perfect forms and Voronoi decompositions to this setting. Our first main result, Theorem~\ref{pf_bounds}, provides an upper bound for the number of $\GLo$-equivalence classes of perfect forms, generalizing van Woerden's result \cite{woerden}. More generally, we obtain uniform bounds for the number of equivalence classes of $k$-dimensional cells of the Voronoi tessellation and for the number of their facets. As an application, we derive in Theorem~\ref{k-bounds} improved bounds for the torsion in $K_m(\mathcal O_F)$, away from a prescribed set of small primes.

We briefly explain how we extend van Woerden's volumetric argument \cite{woerden} to cells of arbitrary dimension. Van Woerden shows that every Voronoi domain of a perfect form admits a representative in its $\GL$-orbit for which all full-dimensional simplices spanned by its vertices have a uniform lower volume bound. Since the trace-one section of the cone is bounded and the interiors of distinct Voronoi domains are disjoint, this yields an upper bound for the number of $\GL$-orbits of perfect forms.

To extend this argument to $k$-dimensional cells, we associate with each such cell a full-dimensional polyhedron constructed through successive pyramids along a maximal flag of faces. Using barycentric subdivisions, these polyhedra can be chosen with pairwise disjoint interiors. Over $F$, instead of choosing a basis of $\rint^n$ consisting of vectors of bounded length, we replace $(\rint^n,Q)$ with an isometric pair $(L,Q^\prime)$, where $L\subseteq\rint^n$ is a finite-index sublattice. This gives a uniform bound for the Euclidean norms of the minimal vectors of $Q^\prime$ on $L$, together with a bound for $[\rint^n:L]$. We then apply the volumetric argument in the tessellation associated with each $L$ and sum over the possible embedded sublattices. 

For the K-theoretic application, our torsion bound takes the following form.
\begin{introtheorem}
Let $m\geq1$, put $n=2m+1$, and let
$$t\coloneqq\operatorname{card}_{b} K_m(\rint)_{\mathrm{tors}}.$$
Then
$$\log(t)\leq e^{O(n^2d(\log n+\log d))}|\disc|^{O(n^2)}.$$
\end{introtheorem}
Here $\operatorname{card}_b(A)$ denotes the order of the subgroup
of $A$ consisting of elements of order prime to all primes at
most $b$.

We now compare our bounds with previous results. Note that we  do not provide the best possible constants, as our objective is to improve the asymptotic order. For $F=\mathbb Q$, our bound improves Soul\'e's estimate by a factor of $n^2$ in the exponent. Moshchevitin's lower bound \cite{moshchevitin2026} shows that the scale $e^{O(n^2\log n)}$ is optimal for the number of perfect-form orbits, and hence for the total cell count used in our application of Gabber's lemma. For totally real and CM fields, the bound of Bayer-Fluckiger, Emery, and Houriet \cite[Theorem~4.3]{Juliet_PHD} involves an exponential factor of order $e^{O(n^4d^4\log n\log d)}$. In addition to improving this factor, our bound has discriminant exponent $O(n^2)$, eliminating the dependence on the degree $d$ in that exponent. For imaginary quadratic fields, our bound recovers the quadratic
dependence on $n$ in the discriminant exponent of Emery's bound
\cite[Proposition~4.1]{Emery_imaginary_quad}, while making the
bound explicit. Our method also applies to totally real and CM
fields with more than one archimedean place, avoiding the
obstruction to extending Emery's approach discussed in
\cite[Remark~4.3]{Emery_imaginary_quad}.
These comparisons are made away from the respective sets of
small primes and are summarized in the following table:

\begin{table}[ht]
\label{table}
\centering
\begin{tabular}{llll}
\hline
Field & Bound & Exponential factor & Discriminant factor \\
\hline
$F=\mathbb Q$ & Soul\'e \cite{soule_bounds} & $e^{O(n^4\log n)}$ & $1$ \\
 & This work & $e^{O(n^2\log n)}$ & $1$ \\
\hline
$F\neq\mathbb Q$ & \small{Bayer et al.} \cite{Juliet_PHD} & $e^{O(n^4d^4\log n\log d)}$ & $|\disc|^{O(n^4d^3\log n)}$ \\
CM or totally real& This work & $e^{O(n^2d(\log n+\log d))}$ & $|\disc|^{O(n^2)}$ \\
\hline
Imaginary quadratic & Emery \cite{Emery_imaginary_quad} & Non-explicit & $|\disc|^{O(n^2)}$ \\
 & This work & Explicit & $|\disc|^{O(n^2)}$ \\
\hline
\end{tabular}
\caption{Factors in the bounds for $\log(\operatorname{card}_{b} K_m(\rint)_{\mathrm{tors}})$ for a suitable $b$. The bounds in \cite{Juliet_PHD} are taken from Houriet's thesis,
but are due to Bayer-Fluckiger, Emery, and Houriet.}
\end{table}
We emphasise the improvement in the case $F=\mathbb Q$, for which the general bound had remained unchanged since Soulé's original work, as noted in \cite{emery, Juliet_PHD}. In addition, for totally real and CM fields, we prove that the discriminant exponent in our bound is independent of the degree. Our estimates nevertheless remain far from the following conjecture of Soulé \cite[Section~5.1]{soule_hermitian_bounds}:
\begin{conjecture}
Fix $d\geq2$. There exist positive constants $\alpha,\beta,\gamma$ such that, for every number field $F$ of degree $d$ over $\mathbb Q$ and every $m\geq1$,
$$|K_m(\rint)_{\mathrm{tors}}|\leq\alpha e^{\beta m^\gamma\log|\disc|}.$$
Moreover, $\gamma$ can be chosen independently of $d$.
\end{conjecture}

The same asymptotic improvement over \cite{Juliet_PHD} is obtained for the number of perfect forms when $F\neq\mathbb Q$. In the rational case, our estimate for perfect forms recovers van Woerden's upper bound.

\section{Preliminaries}

\subsection{Hermitian Matrices over $F_\mathbb{R}$}
Let $F$ be a number field of degree $d$ with ring of integers $\rint$ and discriminant $\disc$. Let $(r_1,r_2)$ be the signature of $F$. Put $d_F(n)=r_1\frac{n(n+1)}{2}+r_2n^2-1$. We write $\fr = F\otimes_\mathbb{Q} \mathbb{R}$ and we let $E_F$ be the set of embeddings of $F$ in $\mathbb{C}.$
The map
$$
F \longrightarrow \mathbb C^{E_F},
\qquad
x\longmapsto (e(x))_{e\in E_F},
$$
induces the Minkowski isomorphism of $\mathbb{R}-$algebras
\begin{equation}
\label{mikowsky}
    F\otimes_{\mathbb Q}\mathbb R
\cong
\mathbb R^{r_1}\times\mathbb C^{r_2}.
\end{equation}

We denote by $c :\fr \to \fr$ the canonical conjugation given by the componentwise conjugation in \eqref{mikowsky}.

Consider 
$$\Her = \{ A\in\mathcal{M}_{n}(\fr) \hspace{5pt} | \hspace{5pt} A^* = A\},$$
where $A^*$ is obtained by applying $c$ entry-wise to $A^t$.
By \eqref{mikowsky},
\begin{equation}
\label{iso_herm_mat}
    \Her \cong \sym^{r_1}\times \mathcal{H}_n(\mathbb{C})^{r_2}.
\end{equation}

We endow $\Her$ with the inner product given by 
$$\langle A,B\rangle = \tr_{F/\mathbb{Q}}(tr(AB)) = \sum_{e\in E_F}tr(A_{e}B_{e}).$$
 
In particular, every $A\in\Her$ defines a quadratic form
$$A(x) = \langle A, xx^*\rangle \in\mathbb{R}$$
for all $x\in \fr^n$. We write $q_A(x)=A(x)$ and $\|x\|^2=\sum_{e\in E_F}x_e^*x_e$. All sums and products over $E_F$ include both embeddings in each complex-conjugate pair.

We equip each of the $r_1$ factors of $\sym$ with the trace inner product
$$
\langle A,B\rangle=\operatorname{tr}(AB),
$$
and each of the $r_2$ factors of $\mathcal H_n(\mathbb C)$ with the inner product
$$
\langle A,B\rangle=2\mathfrak{R}(\operatorname{tr}(AB)).
$$
With respect to these inner products, the isomorphism in \eqref{iso_herm_mat} is an isometry.

Consider the map
$
\phi:\sym
\longrightarrow
\mathbb{R}^{n(n+1)/2}
$
defined by
\begin{equation}
\label{isometry_real}
\phi(h)
=
(h_{11},\sqrt{2}h_{12},\ldots,\sqrt{2}h_{1n},h_{22},\dots).
\end{equation}
and
$
\psi:\mathcal{H}_n(\mathbb{C})
\longrightarrow
\mathbb{R}^{n^2}
$
defined by
\begin{equation}
\label{isometry_complex}
\psi(h)
=
\sqrt{2}\bigl(
h_{11},
\sqrt{2}\mathfrak{R}(h_{12}),
\sqrt{2}\mathfrak{I}(h_{12}),
\ldots,
\sqrt{2}\mathfrak{R}(h_{1n}),
\sqrt{2}\mathfrak{I}(h_{1n}),
h_{22},
\ldots
\bigr),
\end{equation}
where the diagonal entries and the real and imaginary parts of the off-diagonal entries are ordered lexicographically.

Then, the composition of \eqref{iso_herm_mat} with $\phi^{r_1}
\times
\psi^{r_2}$ gives an isometry
\begin{equation}\label{isometry_vor_space_real}
    \Her \cong \mathbb{R}^{r_1\frac{n(n+1)}{2}+r_2n^2}
\end{equation}
with respect to the usual inner product. 

We let $C_{F,n}$ be the cone of positive-definite matrices in $\Her$. Under the above identification
$$C_{F,n} \cong \operatorname{Sym}^+_n(\mathbb{R})^{r_1}\times\mathcal{H}_n^+(\mathbb{C})^{r_2}.$$
We denote by $\overline{C}_{F,n}$ the closure of $C_{F,n}$ in $\Her$, that is, the cone of positive semidefinite forms. We say that the radical of $h\in\overline{C}_{F,n}$ is defined over $F$ if there exists an $F$-subspace $W\subseteq F^n$ such that $Rad(h)=W\otimes_{\mathbb Q}\mathbb R$. The \emph{rational closure} $C_{F,n}^{*}$ of $C_{F,n}$ is the subcone of $\overline{C}_{F,n}$ consisting of those forms whose radical is defined over $F$.

\subsection{Perfect Forms and the Voronoi Tessellation}
We introduce the theory of Perfect forms originally introduced by Voronoi \cite{Voronoi1908} and later generalized to arbitrary number fields by Koecher \cite{Koecher}. As remarked in the introduction some of the properties of the perfect forms in the rational case do not hold in the general one. 

\begin{definition}
\label{def_varias}
Let $h \in C_{F,n}$. We define the following notation:
\begin{enumerate}[label=(\roman*)]
    \item the \textit{minimum} of $h$ by
    $
        \mu(h) := \min_{x \in \rint^n \setminus \{0\}} h(x).
    $
    \item the set of \textit{minimal vectors} of $h$ by
    $
        m(h) := \{ x \in \rint^n \setminus \{0\} \mid h(x) = \mu(h) \}.
    $
    \item For $x\in \fr^n$, we set $\hat{x}=xx^*$ and call it the \textit{Gram matrix} associated with $x$.
    \item We say that $h$ is \textit{perfect} if the $\mathbb{R}-$span of the set
    $
        \widehat{m}(h)=\{ \hat x \mid x \in m(h) \}
    $
    is $\Her$.
    \item We say that $h$ is \textit{well rounded} if $\mu(h) = 1$ and the $F$-span of $m(h)$ is $F^n$.
    \item  We denote by $\perfect{F,n}  \subset C_{F,n}$ the subset of perfect forms with unitary minimum.
\end{enumerate}
\end{definition}

\begin{definition}
Let $h\in C_{F,n}$. We define its \textit{Voronoi domain} as the following set
$$\mathcal{D}(h) = \left\{ \sum_{x\in m(h)} \lambda_x \hat x : \lambda_x \geq 0 \text{ for all } x\in m(h) \right\}.$$
\end{definition}

\begin{remark}
Let $h\in C_{F,n}$. The Voronoi domain $\mathcal{D}(h)$ is a pointed convex polyhedral cone. Each of its extreme rays is generated by an element of the form $xx^*$, for some $x\in m(h)$. By \cite[Theorem~8]{Jaquet_rank_7}, the converse holds when
$F=\mathbb{Q}$, and the same argument applies when $F$ is
an imaginary quadratic field. However, this argument does not extend to general number fields:
when $r_1+r_2>1$, the rays generated by minimal vectors are
not extreme rays of $\overline{C}_{F,n}$, as shown in
Remark~\ref{interior_min_points}.
\end{remark}

\begin{definition}
Let $\sigma \neq \{0\}$ be a face of a Voronoi domain
of a perfect form. We define the \textit{dimension} of $\sigma$ by
$$
\dim(\sigma) := \dim_{\mathbb{R}}
\operatorname{span}_{\mathbb{R}}(\sigma) - 1.
$$
Throughout, dimensions of nonzero Voronoi cones and their
faces are understood according to this convention.

 We denote by 
$\kfacets$ the set of \textit{$k$-codimensional faces of $\sigma$}. 
We write $$\mathcal{F}(\sigma)=\bigcup_{0\leq k\leq \dim \sigma} \kfacets$$ for the set of \textit{faces of $\sigma$}.
We denote 
$$\mathcal{V}_k = \bigcup_{h\in \perfect{F,n}} \mathcal{F}^{d_F(n)-k}(\mathcal{D}(h)).$$

 For a finite-index $\rint$-sublattice $L\subseteq\rint^n$, we use the same definitions with $L$ in place of $\rint^n$, writing $\mathcal{P}_{F,n}(L)$, $m_L(h)$, $\mathcal{D}_L(h)$ and $\mathcal{V}_k(L)$ for the corresponding minimal vectors, Voronoi domains and faces. We put $\mathcal{V}(L)=\bigcup_{k=0}^{d_F(n)}\mathcal{V}_k(L)$. The Voronoi tessellation for $L$ is understood with these conventions.
\end{definition}

Let $\gamma \in \GLo$. We consider the linear action of $\GLo$ on $\Her$ given by
\begin{equation}
\label{action}
    \begin{aligned}
\tilde{\gamma} \colon \Her &\to \Her \\
X&\mapsto \gamma^*  X \gamma
\end{aligned}.
\end{equation}
In particular, $\GLo$ acts over $C_{F,n}$ and $C_{F,n}^*$. Given $h \in \perfect{F, n}$, this action translates at the level of minimal vectors as
$
m(\tilde{\gamma}(h)) = \gamma^{-1} m(h).
$
And at the level of the Voronoi's domain,
\begin{equation}
\label{real_action}
\mathcal{D}(\tilde{\gamma}(h)) = \gamma^{-1} \mathcal{D}(h) (\gamma^{-1})^*= \widetilde{(\gamma^{-1})^*}\left(\mathcal{D}(h)\right).
\end{equation}

We recall the following result from Braun, Coulangeon, Nebe and Schönnenbeck \cite[Theorem~2.4, Corollary 2.5]{Coulangeon}.
\begin{theorem}
\label{vor_th}
The Voronoi domains of perfect forms constitute a locally finite
exact tessellation of the rational closure $C_{F,n}^{*}$. More
precisely:
\begin{enumerate}[label=(\roman*)]
    \item
    $
        C_{F,n}
        \subseteq
        \bigcup_{h\in\perfect n}\mathcal{D}(h)
        =
        C_{F,n}^{*}
        \subseteq
        \overline{C}_{F,n}.
    $

    \item Let $h,h'\in\perfect{F,n}$, If $\mathcal{D}(h)$ and $\mathcal{D}(h')$ share an interior point, then $\mathcal{D}(h) = \mathcal{D}(h')$.

    \item Every facet of $\mathcal{D}(h)$ is a common facet of
    exactly two Voronoi domains of perfect forms.

    \item  Every compact subset of $C_{F,n}$ intersects only finitely many Voronoi domains.

    \item There are only finitely many
    $\GLo$-equivalence classes of perfect forms of
    minimum~$1$.
\end{enumerate}
\end{theorem}

\subsection{The Voronoi-Koecher complex}
In this section, we introduce the Voronoi-Koecher complex, originally defined by Elbaz-Vincent, Gangl and Soulé \cite{2002, philippe_advances} for $F = \mathbb{Q}$ under the name of Voronoi complex and later generalized to imaginary quadratic fields \cite{hermitian} and arbitrary number fields under the name of Voronoi-Koecher complex \cite{ Vor_comp_num_fields}. 

Let $\Gamma = \GLo$. Denote by $\mathcal{C}_{k}$ a set of representatives of the $\Gamma$-orbits of $\mathcal{V}_k.$ Set
$$
\Sigma_{k}^*(\Gamma)= \{ \sigma \in \mathcal{C}_{k} : \operatorname{relint}(\sigma) \cap \partial C_{F,n}^* = \emptyset \}.
$$
For $\sigma \in \Sigma_k^*(\Gamma)$, we denote by $\widehat{m}(\sigma)$ the set of its extreme rays, and by 
$\mathbb{R}(\sigma)$ the vector subspace of $\mathbb{R}^{d_F(n)+1}$ generated by $\widehat{m}(\sigma)$ and we denote by
$$\Gamma_\sigma = \{ \gamma \in \Gamma : \gamma\cdot\sigma = \sigma \}$$
\textit{the stabilizer of $\sigma$}. We define $\Sigma_k(\Gamma) \subset \Sigma_k^*(\Gamma)$ as the subset
consisting of those $\sigma$ such that no element of $\Gamma_\sigma$ inverts the
orientation of $\mathbb{R}(\sigma)$.

Given $\sigma\in\Sigma_{k}^*(\Gamma)$ and $\tau\in\Sigma_{k-1}^*(\Gamma)$ we
let 
$$[\sigma:\tau] = \sum_{\substack{\tau'\in\facets\\
                \tau'\sim_\Gamma\tau}} \eta(\tau,\tau')\varepsilon(\sigma,\tau').$$
where $\eta$ and $\varepsilon$ are signs controlling how the orientations interact with the action of $\Gamma$. See \cite[Section 3.3]{philippe_advances} for more details.

Let $V_k$ be the free abelian group generated by $\Sigma_k$, and let
$d_k : V_k \to V_{k-1}$ be defined on the generators by 
\begin{equation*}
\label{differential}
    d_k(\sigma) = \sum_{\tau\in\Sigma_{k-1}(\Gamma)} [\sigma:\tau]\tau.
\end{equation*}
for all $k\ge n$. It has been shown in \cite{2002, philippe_advances} that $d_k$ is a differential and endows
$(V_*, d_*)$ with a structure of cellular complex, called 
the \textit{Voronoi-Koecher complex of $\Gamma$}, and denoted  $\operatorname{Vor}_\Gamma$.
\subsection{Volume Formulas}
Let $\Delta = [x_0,\dots,x_n]\subset \mathbb{R}^n$ be a simplex. We denote by $\Delta^i$ the facet of $\Delta$ defined by $[x_0,\dots, \hat x_i,\dots,x_n]$. By $$c_{\Delta} = \frac{1}{n+1}\sum_{i = 0}^n x_i$$ the \textit{barycentre} of $\Delta$.

The $n$-dimensional volume of a measurable set $S \subset \mathbb{R}^{n}$, with respect to the standard Lebesgue measure on $\mathbb{R}^{n}$, is denoted by $\operatorname{Vol}_n(S)$.
In particular the $n$-dimensional unit
ball $B_n$ has volume
$
    \operatorname{Vol}_n(B_n)
    =
    \frac{\pi^{n/2}}{\Gamma(n/2+1)},
$
where $\Gamma$ denotes Euler's gamma function. The $n-$simplex $\Delta\subset\mathbb{R}^n$ has volume
$$
    \operatorname{Vol}_n(\Delta)
    =
    \frac{1}{n!}
    \left| \det\bigl((x_i-x_0)_{i \in \{1,\dots,n\}}\bigr) \right|.
$$
More generally, if we consider $\Delta$ embedded in $\mathbb{R}^m$ with $m>n$, we can express the volume of $\Delta$ as follows. Consider the $m\times n$ matrix $W = (x_i^t-x_0^t)_{i\in\{1,\dots,n\}}$. Then
$$\vol_{n}(\Delta) = \frac{1}{n!}\sqrt{\det W^tW}.$$
If $S \subset \mathbb{R}^n$ is a convex measurable set of affine dimension $n-1$ and $p \in \mathbb{R}^n$, we denote by $d^\perp(p,S)$ \textit{the orthogonal distance}, this is the distance from $p$ to the affine hyperplane containing $S$. It holds that
$
    \operatorname{Vol}_n(S * p)
    =
    \frac{d^{\perp}(p, S)}{n} \operatorname{Vol}_{n-1}(S),
$
where $S*p$ denotes the \textit{join} between $S$ and the point $p$.

In particular, 
\begin{equation}
    \label{volumen_k_join}
    \vol_{F,n}(S*p_1*\cdots*p_m) = \vol_{n-m}(S)\prod_{i = 0}^{m-1} \frac{d^{\perp}(p_{i+1}, S*p_1*\cdots*p_i)}{n-i}.
\end{equation}

We call two convex polytopes $S_1,S_2 \subset \mathbb{R}^n$ \textit{essentially
disjoint} if $\operatorname{Int}(S_1) \cap \operatorname{Int}(S_2) = \emptyset$.
In particular for $m$ pairwise essentially disjoint convex polytopes
$S_1,\ldots,S_m \subset \mathbb{R}^n$ it holds that
\begin{equation}
\label{essentially_disjoint}
      \operatorname{Vol}_n\left(\bigcup_{i=1}^N S_i\right)
    =
    \sum_{i=1}^N \operatorname{Vol}_n(S_i).  
\end{equation}

\section{Bounded Minimal Vectors}

\subsection{Normalization over Totally Real and CM Fields}
In what follows of this paper fix a totally real or CM field $F$. Recall that in these cases componentwise conjugation is induced by an involution of $F$ preserving $\rint$. In particular, the adjoint of an integral matrix is again integral. Let
$$\Codiff=\{x\in F:\operatorname{Tr}_\mathbb{Q}^F (x\rint)\subseteq\mathbb Z\}$$
be the \textit{codifferent ideal} and let $\Diff$ be the \textit{different ideal}. For a nonzero integral ideal $I$, write $N(I)=[\rint:I]$ for its norm. The real dual of $(\rint^n,q_Q)$ is isometric to $((\Codiff)^n,q_{Q^{-1}})$. Indeed,
we can consider the isometry $(\rint^n,q_Q)^* \to ((\Codiff)^n,q_{Q^{-1}})$ given by $x\to Qx$.

\begin{lemma}\label{short_different}
There exists $0\neq a\in\Diff$ such that
$$\kappa(a):=\max_{e\in E_F}|e(a)|^2\leq d|\disc|^{4/d}.$$
\begin{proof}
By \cite[Proposition~4.1]{emery}, the covering radius of the ideal lattice $(\Diff,q_0)$, where $q_0(x)=\sum_{e\in E_F}|e(x)|^2$, is at most
$$R=\frac{\sqrt d}{2}|\disc|^{1/d}N(\Diff)^{1/d}=\frac{\sqrt d}{2}|\disc|^{2/d}.$$
A shortest nonzero vector $a\in\Diff$ satisfies $\|a\|\leq2R$. Hence $\kappa(a)\leq\|a\|^2\leq4R^2=d|\disc|^{4/d}$.
\end{proof}
\end{lemma}

Set
$$U_{F,n}=n^3d^3|\disc|^{4/d},\qquad T_{F,n}=n^{nd}d^{nd/2}|\disc|^{2n}.$$
\begin{theorem}\label{Bound_Hermitian_min_vec}
Let $Q\in C_{F,n}$ be well rounded. There exists $E\in\mathcal M_n(\rint)\cap\operatorname{GL}_n(F)$ such that,
\begin{enumerate}[label = \roman*)]
    \item Let $L = E^*\rint^n$ and $Q' = E^{-1}QE^{-*}$. The map $x\mapsto E^*x$ is an isometry from $(\rint^n,q_Q)$ to $(L,q_{Q'})$.
    \item Every $x\in m_L(Q') = E^* m(Q)$ satisfies
\begin{equation}\label{eq:trans-minvec}
\|x\|^2\leq U_{F,n}.
\end{equation}
\item $[\rint^n:L]\leq T_{F,n}$.
\end{enumerate}
\begin{proof}
Choose $a$ as in Lemma~\ref{short_different} and put $\kappa=\kappa(a)$. Banaszczyk’s transference theorem \cite[Theorem~2.1]{transference} gives
$$\lambda_1(\rint^n,q_Q)\lambda_{nd}((\Codiff)^n,q_{Q^{-1}})\leq nd.$$
Since the first factor is $1$, there exist $nd$ real-linearly independent vectors $y_j\in(\Codiff)^n$ such that $q_{Q^{-1}}(y_j)\leq(nd)^2$. Then $z_j=ay_j\in\rint^n$ remain real-linearly independent and satisfy $q_{Q^{-1}}(z_j)\leq\kappa(nd)^2$. They span $F^n$ over $F$, so we may select $n$ of them that are $F$-linearly independent and let $E$ have these columns. Note that in general $z_i$ might not span $\rint^n$.

Since $E^*\in\mathcal{M}_n(\rint)$, then $L=E^*\rint^n\subseteq\rint^n$. Direct substitution gives $q_{Q'}(E^*x)=q_Q(x)$ which proves i).

By construction,
$$Q'^{-1}=E^*Q^{-1}E,\qquad \sum_{e\in E_F}\operatorname{tr}((Q'_e)^{-1})=\sum_{i=1}^nq_{Q^{-1}}(z_i)\leq n^3d^2\kappa.$$
The largest eigenvalue of each $(Q'_e)^{-1}$ is therefore at most $n^3d^2\kappa$. Applying the Rayleigh inequality componentwise and summing gives $\|v\|^2\leq n^3d^2\kappa q_{Q'}(v)$. This proves \eqref{eq:trans-minvec} for the minimal vectors of $(L,q_{Q'})$.

We are going to prove iii) now. Recall that,
$$[\rint^n:L] = |N_{\mathbb{Q}}^F( \det E^*)| = \left|\prod_{e\in E_F} \det E_e^*\right|.$$
Hadamard's inequality and AGM inequality give:
\begin{align*}
    |N_{\mathbb{Q}}^F( \det E^*)|^2\frac{1}{\prod_{e\in E_F} \det Q} &= \prod_{e\in E_F} \det(E_e^*Q^{-1}_eE_e) \leq \prod_{i=1}^n\prod_{e\in E_F}(z_{i,e}^*Q^{-1}_ez_{i,e})\\
    &\leq \prod_{i = 1}^n \left(\frac{\sum_{e\in E_F} z_{i,e}^*Q^{-1}_ez_{i,e}}{d}\right)^d = \prod_{i=1}^n \left(\frac{q_{Q^{-1}}(z_i)}{d}\right)^d\\
    &\leq (\kappa n^2d)^{nd} 
\end{align*}
Thus,
\begin{equation}
\label{bound_norm_E}
    |N_{\mathbb{Q}}^F( \det E^*)|^2\leq (\kappa n^2d)^{nd} \prod_{e\in E_F} \det Q_e. 
\end{equation}
Remains to bound $\prod_{e\in E_F} \det Q_e$. We proceed with an analogous argument as the one above. Choose $n$ $F-$independent minimal vectors $x_i\in m(Q)$ and let $X$ have these columns. Then,
$$\prod_{e\in E_F}\det(X_e^*Q_eX_e)\leq\prod_{i=1}^n\prod_{e\in E_F}(x_{i,e}^*Q_ex_{i,e})\leq\prod_{i=1}^n\left(\frac{\sum_{e\in E_F}x_{i,e}^*Q_ex_{i,e}}{d}\right)^d=d^{-nd}.$$
Since $\det X$ is a nonzero algebraic integer, $|N_{F/\mathbb{Q}} (\det X)|\geq1$, and consequently
\begin{equation}
\label{bound_det_Q}
    \prod_e\det Q_e\leq d^{-nd}.
\end{equation}
Applying \eqref{bound_det_Q} to \eqref{bound_norm_E}, yields
 $$|N_{F/\mathbb{Q}} (\det E)|^2\leq(\kappa n^2)^{nd}.$$
Thus,
$$[\rint^n:L]=|N_{F/\mathbb{Q}} (\det E)|\leq(n\sqrt\kappa)^{nd}\leq T_{F,n}.$$
\end{proof}
\end{theorem}

\begin{remark}
For $F=\mathbb Q$, passage to a finite-index sublattice is unnecessary. By van Woerden's normalization lemma \cite[Lemma~3.3]{woerden}, every positive-definite quadratic form is $\operatorname{GL}_n(\mathbb Z)$-equivalent to a form whose minimal vectors satisfy $\|x\|^2\leq\frac{1}{8}n^3(n+7)$. Thus, in the rational case, we may work throughout with $L=\mathbb Z^n$ and omit the subsequent count over embedded sublattices.
\end{remark}

\begin{lemma}\label{count_embedded_sublattices}
\label{count_embedded_submodules}
For $T\geq1$, the number of finite-index $\rint$-sublattices $L\subseteq\rint^n$ with $[\rint^n:L]\leq T$ is at most $e^dT^{n+1}$.
\begin{proof}

For $s\in\mathbb{C}$ such that $\mathfrak{R}(s)>n$, we define the map
$$\zeta_{\rint^n}(s) = \sum_{L\subset \rint^n} [\rint^n : L]^{-s}$$
where the sum ranges over all the full sublattices of $\rint^n$. 
Let 
\begin{equation}
\label{eq:bound_sub}
    \mathcal{S} = \{L\subseteq\rint^n:[\rint^n:L]\leq T\}.
\end{equation} Then, for $s\in\mathbb{N}$, $s>n$
$$\# \mathcal{S} \leq  \sum_{L\in \mathcal{S}} \left(\frac{T}{[\rint^n:L]}\right)^s \leq T^s\zeta_{\rint^n}(s).$$

For $s>n$ Solomon's identity holds; see \cite{BushnellReiner1980}:
\begin{equation}
\label{solomon_id}
  \zeta_{\rint^n}(s) = \prod_{j=0}^{n-1}\zeta_F(s-j),  
\end{equation}
where $\zeta_F$ denotes the Dedekind zeta function associated to $F$.
Since $\zeta_F(s)\leq\zeta(s)^d$ for $s>1$ and by \eqref{solomon_id}
\begin{equation}
    \zeta_{\rint^n}(n+1) \leq \left(\prod_{j=2}^{n+1}\zeta(j)\right)^d.
\end{equation}
On the other hand, taking logarithms
\begin{equation}
\label{intermediate_solomon}
\sum_{j=2}^{\infty} \log \zeta(j) \leq \sum_{j=2}^{\infty}(\zeta(j)-1)=\sum_{j = 2}^{\infty}\sum_{r=2}^{\infty} r^{-j}=  \sum_{r=2}^{\infty}\frac{1}{r(r-1)}=1
\end{equation}
Thus, combining \eqref{solomon_id} and $\eqref{intermediate_solomon}$, we have
 $$\zeta_{\rint^n}(n+1) \leq e^d.$$
The result follows by applying the preceding inequalities to \eqref{eq:bound_sub} with $s=n+1$.
\end{proof}
\end{lemma}

\section{Bounds on the Geometry of the Voronoi Tessellation}
Fix $n\in\mathbb{N}$. Let $d_F(n)= r_1\frac{n(n+1)}{2}+r_2n^2-1$. Consider $k\in\mathbb{N}$ such that $$dn-1\leq k\leq d_F(n).$$

In what follows, we identify the Voronoi space $\Her$ with its image under the isometry described in \eqref{isometry_vor_space_real} and we work directly in $\mathbb{R}^{d_F(n)+1}$.

\subsection{The Topological Closure of the Cone of Positive-Definite Hermitian Forms}

The following lemma is a straightforward generalization of \cite[Proposition 7.1.4] {euclidean_lattices_martinet}. We reproduce it here for completeness of the text.
\begin{lemma}\label{martinet}
The convex hull of the elements $\hat{x}\in\Her$ with $x\in\fr^n$, which coincides with $\overline{C}_{F,n}$, is contained in the intersection of the circular cone with axis $I_n\in\Her$ and half-angle $\alpha=\arccos\left(\frac{1}{\sqrt{dn}}\right)$ with the half-space $$\{X\in\Her:\langle I_n,X\rangle\geq 0\}.$$
\begin{proof}

The assertion is immediate for $A=0$.
Let $A=(A_e)_{e\in E_F}\in\overline C_{F,n}\setminus\{0\}$,
with eigenvalues $\lambda_1^e,\ldots,\lambda_n^e$. Let $\beta$ be the angle between $A$ and $I_n$. Then
$$
\cos \beta
=
\frac{\langle A,I_n\rangle}
{\|A\|\,\|I_n\|}.
$$

We have the equalities

$$
\langle A,I_n\rangle
=
\sum_e \operatorname{Tr}(A_e)
=
\sum_{e,i}\lambda_i^e,
\qquad
\|I_n\|^2
=
\sum_e\operatorname{Tr}(I_n) = r_1n+2r_2n
=
dn,
$$

and

$$
\|A\|^2
=
\operatorname{Tr}(A^2)
=
\sum_{e,i}(\lambda^e_i)^2.
$$

Consequently,

$$
\cos \beta
=
\frac{\sum_{e,i}\lambda_i^e}
{\sqrt{dn\sum_{e,i}(\lambda_i^e)^2}}
>0.
$$

This shows, on the one hand, that $A$ lies on the side of the
hyperplane $I_n^{\perp}$ that contains the projections,
and, on the other hand, that the inequalities
$\beta\leq\alpha$ and 
$\cos^2\beta\geq\cos^2\alpha$
are equivalent.

Since
$$
\cos^2\beta
=
\frac{1}{dn}
\frac{(\sum_{e,i}\lambda_i^e)^2}
{\sum_{e,i}(\lambda_i^e)^2}
=
\frac{1}{dn}
\left(
1+
2\frac{
\sum_{(i,e)<(j,f)}
\lambda^e_i\lambda^f_{j}
}{
\sum_{i,e}(\lambda^e_i)^2
}
\right).
$$
For $A\neq 0$, equality holds precisely when $A$ has a
rank-one component at a single real place and zero components
at all other archimedean places.
\end{proof}
\end{lemma}
\begin{remark} \label{interior_min_points}
   Note that, for $x\in\rint^n\setminus\{ 0\}$, we have
$$
\operatorname{rad}(\hat{x})=\bigoplus_{e\in E^*_F}\langle x_e\rangle^\perp,
$$
and hence $\operatorname{rank}(\hat{x})=r_1+r_2$. Consequently, the ray generated by $\hat{x}$ is an extreme ray
of $\overline{C}_{F,n}$ if and only if $r_1+r_2=1$,
equivalently, if and only if $F=\mathbb{Q}$ or $F$ is
an imaginary quadratic field.
\end{remark}

\subsection{Volume Bounds for the Voronoi Domains and Their Faces}
Consider the affine hyperplane $T=\langle I_n\rangle^\perp+\frac{1}{dn}I_n$, that is,
$$T=\{Q\in\Her:\sum_{e\in E_F}\operatorname{tr}(Q_e)=1\}.$$
Denote by $B_T\subset T$ the $d_F(n)$-dimensional ball with centre $\frac{1}{dn}I_n$ and radius $\sqrt{\frac{dn-1}{dn}}$. For a finite-index sublattice $L\subseteq\rint^n$, write $\mathcal V_T(L)=\{\tau\cap T:\tau\in\mathcal V(L)\}$. If $\tau\in\mathcal V_k(L)$, then $\tau_t=\tau\cap T$ has dimension $k$. Choose one generating minimal vector $x$ for each extreme ray of $\tau$. The vertices of $\tau_t$ are the corresponding points
$$v(\tau_t)=\left\{\frac{\hat x}{\|x\|^2}\right\}.$$
Indeed, $\langle I_n,\hat x\rangle=\|x\|^2$, so these points lie in $T$.

\begin{remark}\label{integral_Gram_mat}
For $x,y\in L\subseteq\rint^n$, the matrices $\hat x=xx^*$ and $\hat y=yy^*$ are integral Hermitian matrices, and
$$\langle\hat x,\hat y\rangle=\operatorname{Tr}_{F/\mathbb{Q}} \bigl(\operatorname{tr}(xx^*yy^*)\bigr)\in\mathbb Z.$$
Indeed, for every $\alpha\in\rint$, multiplication by $\alpha$ defines a $\mathbb Z$-linear endomorphism
$
\mu_\alpha:\rint\to\rint,
$
whose trace is $\operatorname{Tr}_{F/\mathbb{Q}}(\alpha)$. In particular,
$
\tr_{F/\mathbb{Q}}(\alpha)\in\mathbb Z.
$
\end{remark}
The following lemma adapts van Woerden's  volumetric argument of \cite[Corollary~3.4]{woerden} to faces of arbitrary dimension and to a fixed embedded sublattice.
\begin{lemma}
\label{inscribed_simplex}
Let $L\subseteq\rint^n$ be a finite-index sublattice and let $\sigma$ be a Voronoi domain of a perfect form on $L$ whose minimal vectors satisfy
\begin{equation}\label{norma_vertices}
\|x\|^2\leq U_{F,n}.
\end{equation}
Let $\tau$ be a face of $\sigma$ with $\dim\tau_t=k$. Every $k$-dimensional simplex $\Delta_\tau$ spanned by vertices of $\tau_t$ satisfies
$$\vol_k(\Delta_\tau)\geq l_{F,n,k}\coloneqq\frac{1}{k!U_{F,n}^{k+1}}.$$
\begin{proof}
Consider $\Delta_\tau$ as in the statement and write its vertices as $\frac{\hat x_i}{\|x_i\|^2}$, $0\leq i\leq k$. Put $\Delta=0*\Delta_\tau$. Since the vertices lie in $T$ and are affinely independent, the matrices $\hat x_i$ are linearly independent. Let $W$ have their images under \eqref{isometry_vor_space_real} as columns. By Remark~\ref{integral_Gram_mat}, $W^tW$ is an integral positive-definite matrix, so $\det(W^tW)\geq1$. Consequently,
\begin{equation}\label{eq_1}
\vol_{k+1}(\Delta)=\frac{\sqrt{\det(W^tW)}}{(k+1)!\prod_{i=0}^k\|x_i\|^2}\geq\frac{1}{(k+1)!U_{F,n}^{k+1}}.
\end{equation}
The pyramid formula gives
\begin{equation}\label{eq_2}
\vol_k(\Delta_\tau)=\frac{(k+1)\vol_{k+1}(\Delta)}{d^\perp(0,\Delta_\tau)},
\end{equation}
 Every positive semidefinite $A\in T$ satisfies $\|A\|^2\leq\langle I_n,A\rangle^2=1$, so $d^\perp(0,\Delta_\tau)\leq1$. Combining \eqref{eq_1} and \eqref{eq_2} proves the claim.
\end{proof}
\end{lemma}

We denote by $B_k$ the $k$-dimensional unit ball and put
$$u_{F,n,k}\coloneqq\left(\frac{dn-1}{dn}\right)^{k/2}\vol_k(B_k).$$
\begin{lemma}
\label{bound_volume_face}
Let $0\leq k\leq d_F(n)$ and let $\tau_t$ be a $k$-dimensional face of $\mathcal V_T(L)$. Then
$$\vol_k(\tau_t)\leq u_{F,n,k}.$$
\begin{proof}
By Lemma~\ref{martinet}, $\tau_t$ is contained in $B_T$. Its affine span has dimension $k$, and its intersection with $B_T$ is a ball of radius at most $\sqrt{\frac{dn-1}{dn}}$. The stated volume bound follows.
\end{proof}
\end{lemma}

\subsection{A Bound for the Number of Generators of the Voronoi-Koecher complex and its Differentials through Barycentric Subdivisions}

For each fixed embedded sublattice $L\subset \rint ^n$, we use barycentric subdivisions to associate a full-dimensional polyhedron with a face of any dimension. 

\begin{remark}
\label{real_dimension_hermitian_matrices}
The dimension of $\Her$ satisfies
$$d_F(n)+1=r_1\frac{n(n+1)}2+r_2n^2\leq d\frac{n(n+1)}2.$$
\end{remark}

\begin{lemma}
\label{order_bound}
For $0\leq k\leq d_F(n)$,
$$\frac{u_{F,n,k}}{l_{F,n,k}}\leq e^{O((k+1)(\log n+\log d))}|\disc|^{4(k+1)/d}.$$
\begin{proof}
By definition,
$$\frac{u_{F,n,k}}{l_{F,n,k}}=\frac{k!\pi^{k/2}}{\Gamma(k/2+1)}\left(\frac{dn-1}{dn}\right)^{k/2}U_{F,n}^{k+1}.$$
Since $k\leq d_F(n)=O(n^2d)$ and $U_{F,n}=n^3d^3|\disc|^{4/d}$, the factorial and all remaining factors not involving $\disc$ are at most $e^{O((k+1)(\log n+\log d))}$.
\end{proof}
\end{lemma}

For every $k$-dimensional face $\tau$ of $\mathcal{V}(L)$ for a fixed $L$, choose $\Delta_\tau$ a full-dimensional simplex generated by a subset of $v(\tau_t)$ and denote by $c_\tau$ its barycentre. Write
$$\mathcal{F}_{max}(\tau) = \left\{\tau < \tau_{k+1} < \cdots < \tau_{d_F(n)} \hspace{5pt}|\hspace{5pt} \tau_{k+i}\in \mathcal{V}_{k+i}(L) \text{ for all } i = 1,\dots, d_F(n)-k \right\}$$
the set of \textit{maximal flags} of $\tau$, where $<$ is the partial order given by inclusion. We let $$\mathcal{T}_k = \bigcup_{\tau\in\mathcal{V}_k}\left\{ \tau_t * c_{\tau_{k+1}}* \dots * c_{\tau_{d_F(n)}} \hspace{5pt} | \hspace{5pt} (\tau_{k+1},\dots,\tau_{d_F(n)})\in\mathcal{F}_{max}(\tau) \text{}\right\}$$ be a \textit{barycentre subdivision of dimension $k$} of $\mathcal{V}_T(L)$. It is clear that a barycentre subdivision is formed by essentially disjoints sets. Thus, in particular:
$$\sum_{\tau\in\mathcal{T}_k}\vol_{d_F(n)}(\tau) = \sum_{h\in \perfect {F,n}(L)} \vol_{d_F(n)}(\mathcal{D}_L(h)_t).$$

\begin{lemma}
\label{bound_vol_bar_subd}
Let $\tau<\tau_{k+1}<\cdots<\tau_{d_F(n)}=\sigma$ be a maximal flag in the tessellation for $L$, where the minimal vectors of $\sigma$ satisfy \eqref{norma_vertices}. Then
$$\vol_{d_F(n)}(\tau_t*c_{\tau_{k+1}}*\cdots*c_\sigma)\geq\frac{(k+1)!}{(d_F(n)+1)!}\,l_{F,n,d_F(n)}.$$

\begin{proof}
By \eqref{volumen_k_join}
$$\vol_{d_F(n)}(\tau_t * c_{\tau_{k+1}} * \cdots *  c_{\tau_{d_F(n)}}) = \vol_ k(\tau_t)\prod_{j = k}^{d_F(n)-1} \frac{1}{j+1} d^{\perp}(c_{\tau_{j+1}}, \tau_t * c_{\tau_{k+1}}*\cdots*c_{\tau_{j}}).$$

On the other hand, by definition of flag $\tau_t*c_{\tau_{k+1}}*\cdots*c_{\tau_j}\subset\tau_j$ and since $c_{\tau_{j+1}}\notin\tau_j$, we get
$$d^{\perp}(c_{\tau_{j+1}}, \tau_t * c_{\tau_{k+1}}*\cdots*c_{\tau_{j}}) = d^{\perp}(c_{\tau_{j+1}},(\tau_{j})_t).$$
Let $\Delta_{\tau_j} = [w_{j,0},\dots, w_{j,j}]$ and recall that $c_{\tau_j} = \frac{1}{j+1}\sum_{r=0}^j w_{j,r}$. Since $\tau_{j-1}\in\mathcal{F}^1(\tau_j)$, all the vertices of $\Delta_{\tau_j}$ lie in the same closed half-space bounded by the affine span of $(\tau_{j-1})_t$ within the affine span of $(\tau_j)_t$. Thus
$$d^{\perp}(c_{\tau_j},(\tau_{j-1})_t) = \frac{1}{j+1}\sum_{r=0}^j d^{\perp}(w_{j,r},(\tau_{j-1})_t).$$

Take $v_j$  such that $d^\perp(v_j, (\tau_{j-1})_t)= \max_{r=0,\dots,j} d^\perp (w_{j,r}, (\tau_{j-1})_t)$. Since $\dim \tau_{j} = \dim \tau_{j-1}+1$, it follows that $d^\perp(v_j, (\tau_{j-1})_t)>0$. Combining the previous equalities 
\begin{equation}\label{smart}
    d^{\perp}(c_{\tau_{j+1}}, \tau_t * c_{\tau_{k+1}}*\cdots*c_{\tau_{j}}) \geq \frac{1}{j+2} d^\perp (v_{j+1}, (\tau_j)_t).
\end{equation}

It follows that
\begin{equation}\label{eq1}
  \vol_{d_F(n)}(\tau_t * c_{\tau_{k+1}} * \cdots *  c_{\tau_{d_F(n)}}) \geq  \frac{(k+1)!}{(d_F(n)+1)!} \left(\vol_k(\Delta_{\tau_t}) \prod_{j = k}^{d_F(n)-1} \frac{1}{j+1}d^\perp  (v_{j+1}, (\tau_j)_t) \right).  
\end{equation}

Consider the $d_F(n)$-dimensional simplex $\Delta = \Delta_{\tau_t} * v_{k+1}\dots * v_{d_F(n)}$. By definition, the vertices of $\Delta$ are a subset of $v(\tau_{d_F(n)})$, and thus by Lemma \ref{inscribed_simplex} 
\begin{equation}\label{bounded_vol}
    \vol_{d_F(n)}(\Delta)\geq l_{F,n,d_F(n)}.
\end{equation} 
Then, using \eqref{volumen_k_join} and $v_j\in v(\tau_j)$ we have
\begin{equation}\label{simplex_vertices_min_vec}
    \vol_{d_F(n)}(\Delta) = \vol_k(\Delta_{\tau}) \prod_{j = k}^{d_F(n)-1} \frac{1}{j+1}d^\perp (v_{j+1}, (\tau_j)_t) .
\end{equation}
Applying \eqref{bounded_vol} and \eqref{simplex_vertices_min_vec} to \eqref{eq1} we obtain
$$\vol_{d_F(n)}(\tau_t * c_{\tau_{k+1}} * \cdots *  c_{\tau_{d_F(n)}}) \geq  \frac{(k+1)!}{(d_F(n)+1)!}\cdot l_{F,n,d_F(n)}.$$
\end{proof}
\end{lemma}

Put $\beta_F(n)=2n(n+1)+\frac{4(d_F(n)+1)}{d}$.
\begin{theorem}
\label{bound_number_generators_vor_com}
Let $dn-1\leq k\leq d_F(n)$. Then
$$|\Sigma_k^*|\leq\varepsilon_F(n,k)\coloneqq e^dT_{F,n}^{n+1}\frac{(d_F(n)+1)!}{(k+1)!}\frac{u_{F,n,d_F(n)}}{l_{F,n,d_F(n)}}\leq e^{O(n^2d(\log n+\log d))}|\disc|^{\beta_F(n)}.$$

\begin{proof}
For $\tau\in\Sigma_k^*$, choose $\sigma_\tau$ to be the Voronoi domain of a perfect form containing it. By Theorem~\ref{Bound_Hermitian_min_vec}, there exists $E_\tau\in\mathcal{M}_n(\rint)\cap\operatorname{GL}_n(F)$ such that, setting $L_\tau=E_\tau^*\rint^n$, the image $E_\tau^*\sigma_\tau E_\tau$ is a top cell of the tessellation $\mathcal{V}(L_\tau)$ associated with a perfect form whose minimal vectors satisfy \eqref{norma_vertices}.

Consider $\Sigma_k^*(L)=\{\tau\in\Sigma_k^*:L_\tau=L\}$. Then,
\begin{equation*}
\label{sublatt_sum}
   |\Sigma_k^*|=\sum_{L\subseteq\rint^n}|\Sigma_k^*(L)|.
\end{equation*}
Moreover, distinct elements of $\Sigma_k^*(L)$ have distinct normalized images. Indeed, consider $\tau_1,\tau_2\in\Sigma_k^*(L)$ and denote $E_i=E_{\tau_i}$. Suppose that $(1)$ $E_1^*\rint^n=E_2^*\rint^n$ and $(2)$ $E_1^*\tau_1E_1=E_2^*\tau_2E_2$. Condition $(1)$ implies that $E_2^{-*}E_1^*\in\GLo$. Condition $(2)$ implies that $E_2^{-*}E_1^*\tau_1E_1E_2^{-1}=\tau_2$. Thus, by $(1)$, $\tau_1$ and $\tau_2$ lie in the same $\GLo$-orbit. By definition of $\Sigma_k^*$, this can only happen if $\tau_1=\tau_2$.

Suppose that there exists a uniform upper bound $\Omega$ for the cardinalities of the sets $\Sigma_k^*(L)$. Then, by Theorem~\ref{Bound_Hermitian_min_vec} (iii) and Lemma~\ref{count_embedded_sublattices},
\begin{equation}
    \label{1}
    |\Sigma_k^*|=\sum_{\substack{L\subseteq\rint^n\\ \Sigma_k^*(L)\neq \emptyset} }|\Sigma_k^*(L)| \leq\#\{L\subseteq\rint^n:[\rint^n:L]\leq T_{F,n}\}\Omega\leq e^dT_{F,n}^{n+1}\Omega.
\end{equation}
It remains to find $\Omega$. Fix $L$ and work in the tessellation $\mathcal{V}_T(L)$. By abuse of notation, we continue to denote the normalized images of $\tau$ and $\sigma_\tau$ by $\tau$ and $\sigma_\tau$, respectively, for every $\tau\in\Sigma_k^*(L)$.

Choose the simplices and their barycentres on the faces of the tessellation, and fix the corresponding barycentric subdivision of dimension $k$. For every $\tau\in\Sigma_k^*(L)$, choose a maximal flag $$f_\tau=(\tau,\tau_{k+1},\dots,\tau_{d_F(n)}=\sigma_\tau)\in\mathcal{F}_{max}(\tau)$$ and let $P_\tau=\tau_t*c_{\tau_{k+1}}*\cdots*c_{\sigma_\tau}$ be the associated polyhedron. Every face in this flag belongs to $\sigma_\tau$, so the common bound \eqref{norma_vertices} allows us to apply the volume estimate in Lemma~\ref{bound_vol_bar_subd} to $P_\tau$.

The subdivision pieces are essentially disjoint, and the selected $P_\tau$ are distinct pieces. Hence
$$
\begin{aligned}
\sum_{h\in\mathcal \mathcal{P}_{F,n}(L)} \vol_{d_F(n)}(\mathcal{D}_L(h)_t)
&=
\sum_{\tau \in \mathcal{V}_k(L)}
\sum_{\substack{
(\tau_{k+1},\dots,\tau_{d_F(n)})
\in \mathcal{F}_{max}(\tau)
}}
\vol_{d_F(n)}\!\left(\tau_t * c_{\tau_{k+1}} * \cdots * c_{\tau_{d_F(n)}}\right) \\
&\geq
\sum_{\tau\in\Sigma_k^*(L)}
\vol_{d_F(n)}(P_\tau) \geq |\Sigma_k^*(L)| \cdot \frac{(k+1)!}{(d_F(n)+1)!}\cdot l_{F,n,d_F(n)}.
\end{aligned}
$$

The last inequality follows from Lemma~\ref{bound_vol_bar_subd}.
By Lemma~\ref{martinet}, the sets $\mathcal{D}_L(h)_t$ are contained in $B_T$. Since they are essentially disjoint,
$$\sum_{h\in\mathcal \mathcal{P}_{F,n}(L)} \vol_{d_F(n)}(\mathcal{D}_L(h)_t) \leq \vol_{d_F(n)}(B_T)=u_{F,n,d_F(n)}.$$

We conclude that
$$|\Sigma_k^*(L)|\leq \frac{(d_F(n)+1)!}{(k+1)!}\cdot \frac{u_{F,n,d_F(n)}}{l_{F,n,d_F(n)}}.$$
Thus we can take $\Omega = \frac{(d_F(n)+1)!}{(k+1)!}\cdot \frac{u_{F,n,d_F(n)}}{l_{F,n,d_F(n)}}$ in \eqref{1}, which proves the explicit bound.

The discriminant powers in $T_{F,n}^{n+1}$ and $U_{F,n}^{d_F(n)+1}$ are $2n(n+1)$ and $\frac{4(d_F(n)+1)}{d}$, respectively. All other factors, including the factorials, are at most $e^{O(n^2d(\log n+\log d))}$. This proves the asymptotic estimate.
\end{proof}
\end{theorem}

\begin{theorem}
\label{bound_number_of_facets_vor_dom}
Consider $\sigma\in\Sigma_{k+1}^*$. 
Then, $$|\facets|\leq \delta_F(n,k+1) \coloneqq(k+2)\cdot \frac{u_{F,n,k+1}}{l_{F,n,k+1}} \leq e^{O((k+2)(\log n+\log d))}|\disc|^{\frac{4(k+2)}{d}}.$$
\begin{proof}
Without loss of generality, by Theorem~\ref{Bound_Hermitian_min_vec}, we can work in the tessellation associated with a finite-index sublattice $L\subseteq\rint^n$, where $\sigma$ is identified with a face of the Voronoi domain of a perfect form $\rho$ on $L$ whose minimal vectors satisfy \eqref{norma_vertices}.

Fix $\Delta_\sigma$ a $(k+1)-$dimensional simplex generated by a subset of $v(\sigma_t)$ and let $c_\sigma$ be its barycentre. Then, the set $\{\tau_t*c_\sigma \hspace{5pt} | \hspace{5pt} \tau\in\facets\}$ gives a tessellation of $\sigma_t$, in particular, it is formed by essentially disjoint sets. Thus, applying \eqref{volumen_k_join} and \eqref{essentially_disjoint} we get
$$\sum_{\tau\in\facets} \frac{1}{k+1} \vol_k(\tau_t) \cdot d^{\perp}(c_\sigma,\tau_t) = \sum_{\tau\in\facets} \vol_{k+1}(\tau_t * c_\sigma) = \vol_{k+1}(\sigma_t) \leq u_{F,n,k+1},$$
where the inequality is given by Lemma \ref{bound_volume_face}.

By an analogous argument as in Lemma~\ref{bound_vol_bar_subd}, we get
$$\sum_{\tau\in\facets} \frac{1}{k+1} \vol_k(\tau_t)\cdot d^{\perp}(c_\sigma,\tau_t)\geq \frac{|\facets|\cdot l_{F,n,k+1}}{k+2} ,$$
and thus 
$$|\facets|\leq (k+2)\cdot \frac{u_{F,n,k+1}}{l_{F,n,k+1}}$$
The claim about the order is given by Lemma \ref{order_bound}.
\end{proof}
\end{theorem}
\begin{remark}
    For $F=\mathbb{Q}$, Theorem~\ref{bound_number_of_facets_vor_dom} improves upon Soulé’s bound  $e^{O(kn)}$ \cite{soule_bounds}. The works \cite{emery, Juliet_PHD} use the well-rounded retract directly, so this approximation was not performed.
\end{remark}

\section{Two Consequences}
\subsection{A Bound for the Number of Hermitian Perfect Forms}
\begin{theorem}
\label{pf_bounds}
Let $n\geq2$. The number $p_{F,n}$ of $\GLo$-equivalence classes of perfect Hermitian forms of minimum $1$ over $F$ satisfies
$$p_{F,n}\leq e^dT_{F,n}^{n+1}\frac{u_{F,n,d_F(n)}}{l_{F,n,d_F(n)}}\leq e^{O(n^2d(\log n+\log d))}|\disc|^{\beta_F(n)}.$$
\begin{proof}
Apply the counting argument in Theorem~\ref{bound_number_generators_vor_com} to top-dimensional domains, taking $k=d_F(n)$. 
\end{proof}
\end{theorem}

\begin{remark}
For $F=\mathbb Q$, Theorem~\ref{pf_bounds} recovers the asymptotic order of van Woerden's bound \cite[Theorem~3.1]{woerden}. For $F\neq\mathbb Q$, the bound obtained by
Bayer-Fluckiger, Emery, and Houriet using
the well-rounded retract is
$$p_{F,n}\leq e^{O(n^4d^4\log n\log d)}|\disc|^{O(n^4d^3\log n)};$$
see \cite[Lemma~2.15 and Lemmas~4.1 and~4.2]{Juliet_PHD}. Note that our method gives a discriminant exponent $\beta_F(n)\leq4n(n+1)$, of order $n^2$ independently of the degree.
\end{remark}

\subsection{A Bound for Torsion Part of the K-theory of $\rint$}
For the convenience of the reader, we recall some facts from Soulé \cite{soule_bounds, soule_hermitian_bounds} about the relationship between the Voronoi-Koecher complex associated to $\GLo$ and the algebraic K-theory of $\rint$.

\begin{lemma}
\label{improved_cutoff}
Every prime $p$ for which $\GLo$ contains an element of
order $p$ satisfies $p\le dn+1$.
If, moreover, $F\cap\mathbb Q(\zeta_p)\in
\{\mathbb Q,\mathbb Q(\zeta_p)\}$
for every prime $p$, then every such prime satisfies
$p\le1+\max\{d,n\}$.
\begin{proof}
    Choosing a $\mathbb Z$-basis of $\rint$, we can identify $\GLo$ with a subgroup of $\operatorname{GL}_{nd}(\mathbb Z)$. Thus, if $p$ is a prime such that there exists an element $g\in\GLo$ of order $p$, we have $p\leq dn+1$; see, for instance, \cite[Lemma~7.1]{philippe_advances}.

    Suppose now that $F\cap\mathbb Q(\zeta_p)\in
\{\mathbb Q,\mathbb Q(\zeta_p)\}$ for every prime $p$.
Consider $I_n\neq\gamma\in\GLo$ such that
$\gamma^p=I_n$. Equivalently, the minimal polynomial
$m_\gamma$ of $\gamma$ divides $x^p-1$.
Since $x^p-1$ has distinct roots, so does $m_\gamma$.
Thus $\gamma$ is diagonalizable over an algebraic closure
of $F$. As $\gamma\neq I_n$, it has an eigenvalue
$\alpha\neq1$. Since $\alpha^p=1$ and $p$ is prime,
$\alpha$ is a primitive $p$-th root of unity.
Therefore $P(\alpha)=0$, where
$$
P(x)=x^{p-1}+x^{p-2}+\dots+1
$$
is the $p$-th cyclotomic polynomial.
If $F\cap\mathbb Q(\zeta_p)=\mathbb Q$, then $P$ is
irreducible over $F$, and therefore it divides the
characteristic polynomial of $\gamma$. Thus $p-1\le n$.
Otherwise, our hypothesis gives $\mathbb Q(\zeta_p)\subseteq F$,
and hence $p-1\le d$.
\end{proof}
\end{lemma}

\begin{remark}
Lemma~\ref{improved_cutoff} corrects an imprecision in
\cite[Lemma~9]{soule_hermitian_bounds}, whose proof assumes
that $P$ is irreducible over $F$ whenever
$\mathbb Q(\zeta_p)\not\subseteq F$. This need not hold
when $F\cap\mathbb Q(\zeta_p)$ is a nontrivial proper
subfield of $\mathbb Q(\zeta_p)$.
The same prime cutoff is used in \cite{emery,Juliet_PHD}.
This imprecision affects only the prime cutoff; the remaining
arguments in these works are unchanged.
\end{remark}

Fix $m\geq1$, put $n=2m+1$, and let
$\Gamma=\GLo$. Fix $b=dn+1$, or $b=\max\{d,n\}+1$ under the hypothesis of Lemma~\ref{improved_cutoff}.
Let $\mathcal{S}_b$ be the Serre class of finite abelian
groups whose orders are divisible only by primes at most
$b$. For a finite abelian group $A$, write
$\operatorname{card}_b(A)$ for the order of its subgroup
consisting of elements of order prime to all primes at
most $b$.

Since $n$ is odd, the action of $\GLo$ on the real vector space $\Her$ preserves
orientation. Indeed, on a real factor, the orientation
character is trivial \cite[Lemma~7.2]{philippe_advances},
and the complex factors also preserve orientation \cite[Remark~2.5]{hermitian}.
Thus the orientation module in
\cite[Theorem~2.2]{Vor_comp_num_fields} is trivial.
Consequently, modulo $\mathcal{S}_b$, we have
\begin{equation}
    \label{comparasion}
    H^k(\GLo, \mathbb{Z}) \cong H_{d_F(n)-k}(\operatorname{Vor}_{\GLo},\mathbb{Z}).
\end{equation}
For any $m\geq1$, we let $K_m(\rint)$ be the $m$-th algebraic K-group of $\rint$. We consider the Hurewicz map
$$K_m(\rint)\to H_m(\operatorname{GL}(\rint),\mathbb Z).$$
By Arlettaz’s theorem, its kernel has no $p$-torsion for primes
$p>\frac{m+1}{2}$; see \cite[Theorem~1.5]{Hurewicz}.
Moreover, homological stability gives
$$H_m(\operatorname{GL}(\rint),\mathbb Z)\cong H_m(\GLo,\mathbb Z)$$
since $n=2m+1$; see \cite[Theorem~4.11]{stability}.

Fix $k=d_F(n)-m$. Equip the chain groups with the Euclidean norms for which their cellular bases are orthonormal. Each boundary coefficient is a sum of facet incidence signs, so
$$\|d_k(\sigma)\|_2\leq\sum_{\tau\in\Sigma_{k-1}}|[\sigma:\tau]|\leq|\mathcal F^1(\sigma)|\leq\delta_F(n,k).$$
Gabber's lemma \cite{soule_bounds} and Theorem~\ref{bound_number_generators_vor_com} therefore give
\begin{equation}\label{eq_vor_homology_bound}
|H_{k-1}(\operatorname{Vor}_{\GLo},\mathbb{Z})_{\mathrm{tors}}|\leq¡¡ \delta_F(n,k)^{\varepsilon_F(n,k)}.
\end{equation}

By the universal coefficient theorem, the torsion
subgroups of $H_m(\Gamma,\mathbb{Z})$ and
$H^{m+1}(\Gamma,\mathbb{Z})$ have the same order.
Applying \eqref{comparasion} in degree $m+1$ and using the Hurewicz map, whose kernel has no $p$-torsion for $p>b$, therefore gives
\begin{equation}
\label{eq_k_th_bound}
\operatorname{card}_b K_m(\rint)_{\mathrm{tors}}
\leq \delta_F(n,k)^{\varepsilon_F(n,k)}.
\end{equation}

\begin{theorem}
\label{k-bounds}
Let $m\geq1$, put $n=2m+1$, and let
$$t\coloneqq\operatorname{card}_{b} K_m(\rint)_{\mathrm{tors}}.$$
Then
$$\log(t)\leq e^{O(n^2d(\log n+\log d))}|\disc|^{O(n^2)}.$$
\begin{proof}
Put $k=d_F(n)-m$. By \eqref{eq_k_th_bound},
$$\log(t)\leq\varepsilon_F(n,k)\log\delta_F(n,k).$$
Theorem~\ref{bound_number_generators_vor_com} bounds $\varepsilon_F(n,k)$ by $e^{O(n^2d(\log n+\log d))}|\disc|^{\beta_F(n)}$. Moreover, by Theorem~\ref{bound_number_of_facets_vor_dom},
$$\log\delta_F(n,k)=O\bigl(n^2d(\log n+\log d)+n^2\log|\disc|\bigr).$$
Absorbing the polynomial factors into the exponential gives
$$\log(t)\leq e^{O(n^2d(\log n+\log d))}|\disc|^{\beta_F(n)}(1+\log|\disc|).$$
Since $\beta_F(n)\leq4n(n+1)$ and $1+\log|\disc|\leq2|\disc|$, the discriminant factor is $|\disc|^{O(n^2)}$, with the constant factor absorbed into the exponential.
\end{proof}
\end{theorem}

\begin{remark}
For $F=\mathbb Q$, Theorem~\ref{k-bounds} improves Soulé's asymptotic bound
\cite{soule_bounds} for the logarithm of the torsion order from
$e^{O(n^4\log n)}$ to $e^{O(n^2\log n)}$. For totally real and CM fields with $F\neq\mathbb Q$,
the discriminant exponent in the bound of
Bayer-Fluckiger, Emery, and Houriet
\cite[Theorem~4.3]{Juliet_PHD} is 
$O(n^4d^3\log n)$, whereas ours is $O(n^2)$.
The exponential factor is reduced from $e^{O(n^4d^4\log n\log d)}$
to $e^{O(n^2d(\log n+\log d))}$.

For imaginary quadratic fields, our bound has a discriminant exponent
of order $n^2$, recovering the quadratic dependence on $n$ in the
discriminant exponent of Emery's bound
\cite[Proposition~4.1]{Emery_imaginary_quad}, while providing an explicit
bound. These comparisons are made away from the respective sets of
small primes.
\end{remark}
\section*{Acknowledgments}

\thanks{\small The author is fully supported by the COGENT project which has received funding from the European Union's Horizon Europe Programme under the Marie Sklodowska-Curie actions HORIZON-MSCA-2023-DN-01 call (Grant agreement ID: 101169527), and from UK Research and Innovation. The author thanks Philippe Elbaz-Vincent for his guidance throughout the development of this work.

\section*{Usage of LLM Tools}
After the original mathematical arguments had been completed, ChatGPT was used to improve the English of the paper and discuss the bounds obtained by the method used. This led the model to suggest Solomon's identity, which was unknown to the author, to bound the number of finite-index submodules of $\rint^n$ in Lemma~\ref{count_embedded_submodules}. Incorporating this identity into the existing argument allowed the author to avoid the use of a bounded basis argument. This led to a further reduction by a factor of $d\log n$ in the exponent of $|\disc|$ and by a factor of $\log n$ in the exponent of the exponential factor in the bound for the logarithm of the torsion order in algebraic K-theory for $F\neq\mathbb{Q}$. This suggestion did not improve the bound in the case $F=\mathbb{Q}$. 

\newpage

\bibliographystyle{plain}  
\bibliography{referencias}
\end{document}